\documentclass[a4paper,12pt]{article}
\usepackage{amsmath,amssymb,amsbsy,amsfonts,amsthm,latexsym,amsopn,amstext,
            amsxtra,euscript,amscd,verbatim,epsf,colordvi,graphics}
\usepackage{graphicx}

\newtheorem{thm}{Theorem}[section]
   
\newtheorem{cor}[thm]{Corollary} 
\newtheorem{prop}[thm]{Proposition}

\newtheorem{defn}[thm]{Definition}
\newtheorem{conj}[thm]{Conjecture}

\newtheorem{opb}[thm]{Open Problem}

\def \C{\mathbb{C}}

\def \N{\mathbb{N}}

\def \R{\mathbb{R}}

\def \Z{\mathbb{Z}}

\def\cA{{\mathcal A}}

\def\cC{{\mathcal C}}

\def\cI{{\mathcal I}}
\def\cJ{{\mathcal J}}

\def\cO{{\mathcal O}}
\def\cP{{\mathcal P}}

\def\cR{{\mathcal R}}

\def\cX{{\mathcal X}}
\def\cY{{\mathcal Y}}
\def\cZ{{\mathcal Z}}

\def\fl#1{\left\lfloor#1\right\rfloor}
\def\rf#1{\left\lceil#1\right\rceil}

\def\({\left(}
\def\){\right)}
\def\[{\left[}
\def\]{\right]}
\def\<{\langle}
\def\>{\rangle}

\begin{document}

\title{\bf On Whitney numbers of the Order Ideals of Generalized Fences and Crowns} 

\author{
{\sc Alessandro Conflitti}
\thanks{Fully supported by European Research Training Network
``Algebraic Combinatorics in Europe'' under the grant $\#$
HPRN--CT--2001--00272.} \\
{Fakult\"{a}t f\"{u}r Mathematik} \\
{Universit\"{a}t Wien} \\
{Nordbergstra{\ss}e 15} \\
{A-1090 Wien} \\
{Austria} \\
{\tt alessandro.conflitti@univie.ac.at}
}

\date{}

\pagenumbering{arabic}
\maketitle

{\bf AMS Subject Classification:} 05A15, 06A07. 

\begin{abstract}
We solve some recurrences given by E. Munarini and N. Zagaglia Salvi 
proving explicit closed formulas for Whitney numbers of the distributive 
lattices of order ideals of the fence poset and crown poset.
Moreover, we get explicit closed formulas for Whitney numbers of lattices 
of order ideals of fences with higher asymmetric peaks.
\end{abstract}

\section{Introduction and Preliminaries}

In~\cite{MZ} authors consider the distributive lattices of all order ideals 
of the fence poset and crown poset ordered by inclusion, and they are able to prove 
recursive formulas for their Whitney numbers. 
In this paper, using purely combinatorial methods, we solve these recursions giving 
explicit closed formulas for the corresponding rank polynomials. 
Moreover, in \S~\ref{S:as} we consider a more general class of fence posets, namely 
fences with higher asymmetric peaks, and we get explicit closed formulas for Whitney 
numbers of lattices of their order ideals.

For others combinatorial results about lattices of order ideals of finite posets 
and their Whitney numbers, we remind to~\cite{CGW,Gri,Ste1,Ste2}.

In the sequel we collect some definitions, notations and results that will 
be used in the following. 
For $x \in \R$ we let $\fl{x}=\max\{n \in \Z : n \leq x\}$ and 
$\rf{x}=\min\{n \in \Z : n \geq x\}$; for any $n,m \in \N$, $n \leq m$, we 
let $\[n,m\]=\{t \in \N : n \leq t \leq m\}$, and $\[n\]=\[1,n\]$, therefore 
$\[0\]= \emptyset$. For any complex number $a$, we define the \emph{rising 
factorial} as $\(a\)_{0}=1$ and $\(a\)_{m}=\prod_{j=0}^{m-1}\(a+j\)$ for any 
$m \in \N \setminus \{0\}$.The cardinality of a set $\cX$ will be denoted by 
$\# \cX$. 

We follow~\cite{And,GrKl,Sta2} for combinatorics notations and 
terminology. 
We recall that a \emph{ranked poset} is a poset $P$ with a function 
$\rho:P \longrightarrow \N$, called rank, such that $\rho\(y\)=\rho\(z\)+1$ 
whenever $z$ is covered by $y$ in $P$ and $\min\{\rho\(z\) : z \in P\}=0$. 
The \emph{rank polynomial} of a ranked finite poset $P$ is the polynomial 
$$\sum_{z \in P}X^{\rho\(z\)}=\sum_{j \geq 0} \omega_{j} X^{j},$$ 
where $\omega_{j}=\# \{z \in P : \rho\(z\)=j\}$ are called 
\emph{Whitney numbers} of $P$.

An \emph{order ideal} of a poset $P$ is a subset $I \subset P$ such that if 
$y \in I$ and $z \leq y$, then $z \in I$; it is well known that the set of 
all order ideals of $P$ ordered by inclusion is closed under unions and 
intersections, and hence forms a distributive lattice: we denote it by 
$\cJ\(P\)$, viz. $\cJ\(P\)=\{I \subset P : I \text{ is an order ideal}\}$. 
It is not hard to see that its rank function is the cardinality of order 
ideals. 

Given a finite poset $\(P,\leq\)$, we denote with $W_{P}\(k\)$ the 
$k$--th Whitney numbers of the ranked poset of all order ideals of $P$, i.e.
$W_{P}\(k\)=\# \{I \in \cJ\(P\) : \rho\(I\)=j\}$, where $\rho$ is the rank
function of $\cJ\(P\)$, and the rank polynomial of $\cJ\(P\)$ is denoted 
by $\cR_{P}\(X\)$, i. e. $\cR_{P}\(X\)=\sum_{k \geq 0} W_{P}\(k\) X^{k}$.

We denote by $\cZ_{n}$ the \emph{fence} poset of order $n$, viz. the poset 
$\{z_{1}, \ldots ,z_{n}\}$ in which 
$z_{2j-1} \triangleleft z_{2j} \triangleright z_{2j+1}$, for all 
$j \geq 1$, are the cover relations, by $\cI_{n}\(k\)$ the set of order 
ideals of $\cZ_{n}$ with cardinality $k$, and by $f_{n,k}$ the Whitney 
numbers of the poset of all order ideals of a fence of order $n$, viz. 
$f_{n,k} = \# \cI_{n}\(k\)=W_{\cZ_{n}}\(k\)$. \\
We denote by $\cY_{n}$ the \emph{crown} poset of order $2n$, viz. the poset 
$\{\zeta_{0}, \ldots ,\zeta_{2n-1}\}$ in which the cover relations are 
the following: for any $h \in \{0, \ldots n-1\}$ and $k \in \[n\]$, 
$\zeta_{2h} \triangleleft \zeta_{2k-1}$ if and only if 
$\left|2h-2k+1\right| \equiv 1 \pmod{2n}$, therefore 
$$\rho\(\zeta_{j}\)=\begin{cases}
0 & \text{ if } j \equiv 0 \pmod{2}, \\
1 & \text{ if } j \equiv 1 \pmod{2}.
\end{cases}$$
We also denote by $\cO_{n}\(k\)$ the set of order ideals of $\cY_{n}$ with 
cardinality $k$, and by $c_{n,k}$ the Whitney numbers of the poset of all 
order ideals of a crown of order $2n$, viz. 
$c_{n,k} = \# \cO_{n}\(k\)=W_{\cY_{n}}\(k\)$.

Finally we recall, gluing together, Propositions $1$, $3$ and $5$ of~\cite{MZ}, 
which give recursions for the sequences $f_{n,k}$ and $c_{n,k}$.

\begin{prop} \label{P:5}
For any integer $n$ the recurrence identity
$$\begin{cases}
f_{2n,k} & = f_{2n-1,k} + f_{2n-2,k-2} \\
f_{2n+1,k} & = f_{2n,k-1} + f_{2n-1,k} 
\end{cases}$$
holds, where
$$\begin{cases}
f_{n,k} = 0 & \qquad \text{ if } \, k \notin \[0,n\] \; \text{ or } \, 
                       n < 0 \\ 
f_{n,0} = 1 & \qquad \text{ for all } \, n \in \N 
\end{cases}$$
are the initial values.

Moreover, with the same initial values the formula 
$$f_{n+4,k+2}=f_{n+2,k+2}+f_{n+2,k+1}+f_{n+2,k}-f_{n,k}$$
holds, for all $0 \leq k \leq n \in \N$.

Furthermore, 
\begin{eqnarray*}
c_{n+2,k+3} & = & f_{2n+3,k+3} + f_{2n+1,2n+1-k} \\
c_{n+2,k+2} & = & c_{2n+1,k} + f_{2n+3,k+2} - f_{2n-1,k} \\
c_{n+2,k+2} & = & f_{2n+4,k+2} - f_{2n,k}
\end{eqnarray*}
hold, for all $n \in \N$ and all $0 \leq k \leq 2n$.
\end{prop}

\section{Closed Formulas for Whitney Numbers} 

\label{S:cf}

We need the following Proposition, whose proof can be found 
in~\cite{Sta2}. 

\begin{prop} \label{P:s} 
For all non--negative integers $k \leq n$,
\begin{eqnarray*}
\# \{\mathbf{\overline{x}}=\(x_{1}, \ldots ,x_{k}\) \in
   \(\N \setminus \{0\}\)^{k} : \sum_{j=1}^{k} x_{j} = n\} & = &
   \binom{n-1}{k-1}, \\
\# \{\mathbf{\overline{x}}=\(x_{1}, \ldots ,x_{k}\) \in \N^{k} : 
   \sum_{j=1}^{k} x_{j} = n\} & = & \binom{n+k-1}{k-1}, \\ 
\# \{\mathbf{\overline{x}}=\(x_{1}, \ldots ,x_{k}\) \in \N^{k} : 
   \sum_{j=1}^{k} x_{j} \leq n\} & = & \binom{n+k}{k},    
\end{eqnarray*}
hold.
\end{prop}

\begin{thm} \label{T:main}
For all $k,v \in \N$ such that $k \leq 2v+1$, 
$$f_{2v+1,k} = \binom{v+1}{k} + \sum_{j \geq 1} 
\frac{\(k-2j+1\)_{j-1}\(v+j-k+2\)_{k-2j}}{j!\(k-2j-1\)!}$$
holds.
\end{thm}
\begin{proof}
For all integers $0 \leq k \leq n$, we can write 
$$f_{n,k} = \# \cI_{n}\(k\) = \sum_{j \geq 0} \cA\(n,k,j\) = 
\sum_{j=0}^{\min\{\fl{\frac{n}{2}},\rf{\frac{k}{2}}-1\}} \cA\(n,k,j\),$$
where 
$$\cA\(n,k,j\)=\# \{\cJ \in \cI_{n}\(k\) : 
\# \{x \in \cJ : \rho\(x\)=1\}=j\};$$ 
thus we have that $f_{n,0}=1$,
$f_{n,1} = \# \{x \in \cZ_{n} : \rho\(x\)=0\} = \fl{\frac{n+1}{2}} = 
\rf{\frac{n}{2}}$, and $\cA\(n,k,0\) = \binom{f_{n,1}}{k}$. 

Consider a fence $\cZ_{n}$ with odd cardinality, i.e. $n=2v+1$
for some $v \in \N$, and write it as the poset 
$\{z_{1}, \ldots ,z_{2v+1}\}$ in which 
$z_{2\alpha-1} \triangleleft z_{2\alpha} \triangleright z_{2\alpha+1}$, 
for all $\alpha \geq 1$, are the cover relations. \\
For any given $\cJ \in \cI_{2v+1}\(k\)$ (with $k \geq 2$) such that 
$\# \{x \in \cJ : \rho\(x\)=1\}=j \geq 1$ we can split the set 
$\{x \in \cJ : \rho\(x\)=1\}$ in $r$ separated non-empty blocks 
$\cX_{1}, \ldots ,\cX_{r}$, such that $\sum_{t=1}^{r} \# \cX_{t} =j$; 
$z_{2a},z_{2b} \in \cJ$ with $1 \leq a < b \leq v$ are in the same block 
if and only if $z_{2c} \in \cJ$ for all $c$ such that $1 \leq a < c < b \leq v$.
Each $\cX_{t}$ determines $2 \# \cX_{t} +1$ elements in $\cJ$, so this 
decomposition fix $\sum_{t=1}^{r} \(2 \# \cX_{t} +1\) = 2j+r$ 
elements of $\cJ$ ($j$ of these have rank $1$, and the others $j+r$ have rank 
$0$), and obviously the others can be chosen in 
$\binom{v+1-\(j+r\)}{k-\(2j+r\)}$ ways between the remainder elements with 
rank $0$. \\ 
Moreover, the number of such decompositions $\cX_{1}, \ldots ,\cX_{r}$ is 
$\# \cC\(j,r\)$ times the the total numbers of shifts of all blocks 
$\cX_{1}, \ldots ,\cX_{r}$, which can be evaluated in the following way: 
at least one element of rank $1$ has to be into the slot between the blocks 
$\cX_{t}$ and $\cX_{t+1}$, for any $t \in \[r-1\]$, and the others
$v-\(j+r-1\)$ elements can be freely distributed into the $r+1$ slots, viz. 
before $\cX_{1}$, between $\cX_{t}$ and $\cX_{t+1}$, for any 
$t \in \[r-1\]$, and after $\cX_{r}$, thus from Proposition~\ref{P:s} 
$\binom{v-j+1}{r}$ is the searched value.

\noindent
Therefore if we define $\cC\(\mu,\nu\)=
\{\mathbf{\overline{x}}=\(x_{1}, \ldots ,x_{\nu}\) \in
   \(\N \setminus \{0\}\)^{\nu} : \sum_{j=1}^{\nu} x_{j} = \mu\}$
for any $1 \leq \nu \leq \mu$, from Proposition~\ref{P:s} we have
$\# \cC\(\mu,\nu\) = \binom{\mu-1}{\nu-1}$, hence for any $j \geq 1$ 
\begin{eqnarray*}
\cA\(2v+1,k,j\) & = & \sum_{r=1}^{j} 
\sum_{\substack{\mathbf{\overline{x}} \in \cC\(j,r\) \\ 2j+r \leq k \\
j+r-1 \leq v}} \binom{v-j+1}{r} \binom{v+1-\(j+r\)}{k-\(2j+r\)} \\
& = & \sum_{r=1}^{j} 
\sum_{\mathbf{\overline{x}} \in \cC\(j,r\)} \binom{v-j+1}{r} 
\binom{v+1-\(j+r\)}{k-\(2j+r\)} \\
& = & \sum_{r=1}^{j} \binom{j-1}{r-1} \binom{v-j+1}{r} 
\binom{v+1-\(j+r\)}{k-\(2j+r\)}.
\end{eqnarray*}

Therefore we have 
\begin{eqnarray*} 
\lefteqn{f_{2v+1,k} = \binom{v+1}{k} + 
\sum_{j=1}^{\min\{v,\rf{\frac{k}{2}}-1\}} \sum_{r=1}^{j} \binom{j-1}{r-1} 
\binom{v-j+1}{r} \binom{v+1-\(j+r\)}{k-\(2j+r\)}} \\ 
& & = \binom{v+1}{k} + \sum_{j \geq 1} 
\sum_{r \geq 1} \binom{j-1}{j-r} \binom{v-j+1}{r} 
\binom{v+1-\(j+r\)}{k-\(2j+r\)} 
\end{eqnarray*} 

Writing the sum over $r$ in hypergeometric notation and applying 
Chu--Vandermonde summation, see~\cite{Fin,GR,Krat}, we get 
\begin{flalign*}
& \sum_{r \geq 1} \binom{j-1}{j-r} \binom{v-j+1}{r} 
\binom{v+1-\(j+r\)}{k-\(2j+r\)} \\
& = \frac{{}_{2}F_{1}\[\begin{matrix} 1-j,1+2j-k\\ 2 \end{matrix};1\]
\(v+j-k+2\)_{k-2j}}{\(k-2j-1\)!} \\
& = \frac{\(k-2j+1\)_{j-1}\(v+j-k+2\)_{k-2j}}{j!\(k-2j-1\)!},
\end{flalign*}
and the desired result follows. 
\end{proof} 

\begin{cor}
For any $v \in \N$ and all $0 \leq k \leq 2v+1$, the sequence 
$f_{2v+1,k}$ is increasing in $v$, viz. $f_{2\(v+1\)+1,k} > f_{2v+1,k}$. 
\end{cor}
\qed

\begin{defn} \label{D:star}
Let $\(P_{1}, \leq_{1}\)$, $\(P_{2}, \leq_{2}\)$ be finite posets with
cover relations $\triangleleft_{1}$ and $\triangleleft_{2}$, respectively, 
and let $x_{1} \in P_{1}$, $x_{2} \in P_{2}$ be minimal elements. \\
We consider a new element $\widetilde{x}$ which does not belong
to $P_{1} \biguplus P_{2}$ and we define a new poset 
$\(P_{1}\(x_{1}\) \circledast P_{2}\(x_{2}\), \leq\)$ with 
cover relations $\triangleleft$, where 
$$P_{1}\(x_{1}\) \circledast P_{2}\(x_{2}\)=
P_{1} \biguplus P_{2} \biguplus \{\widetilde{x}\},$$ 
and for any $x,y \in P_{1}\(x_{1}\) \circledast P_{2}\(x_{2}\)$ we have 
$x \triangleleft y$ if and only if one of the following conditions holds: 
\begin{description}
\item{$\circ$} $x,y \in P_{1}$ and $x \triangleleft_{1} y$ in $P_{1}$,
\item{$\circ$} $x,y \in P_{2}$ and $x \triangleleft_{2} y$ in $P_{2}$,
\item{$\circ$} $x_{1} \triangleleft \widetilde{x}$,
\item{$\circ$} $x_{2} \triangleleft \widetilde{x}$.
\end{description}
\end{defn}

\begin{thm} \label{T:1fl}
Let $\(P_{1}, \leq_{1}\)$, $\(P_{2}, \leq_{2}\)$ be finite posets, 
$x_{1} \in P_{1}$, $x_{2} \in P_{2}$ be minimal elements, and
$\widehat{P}=P_{1}\(x_{1}\) \circledast P_{2}\(x_{2}\)$; 
then 
$$\cR_{\widehat{P}}\(X\) = \cR_{P_{1}}\(X\) \cR_{P_{2}}\(X\) 
 + X^{3} \cR_{P_{1} \setminus \{x_{1}\}}\(X\) 
\cR_{P_{2} \setminus \{x_{2}\}}\(X\)$$
holds. 
\end{thm}
\begin{proof}
Let us write $\cJ\(\widehat{P}\)= \biguplus_{k=0}^{\# \widehat{P}}
\cJ_{k}$, where $\cJ_{k}=\{I \in \cJ\(\widehat{P}\) : \rho\(I\)=\# I =k\}$,
thus
$$W_{\widehat{P}}\(k\)=\# \cJ_{k}=
\# \{I \in \cJ_{k} : \widetilde{x} \notin I\} +
\# \{I \in \cJ_{k} : \widetilde{x} \in I\}.$$
It is not hard to see that
$$\# \{I \in \cJ_{k} : \widetilde{x} \notin I\}=
\sum_{j=0}^{k} W_{P_{1}}\(j\) \cdot W_{P_{2}}\(k-j\)$$
and 
\begin{eqnarray*}
\# \{I \in \cJ_{k} : \widetilde{x} \in I\} & = &
\# \{I \in \cJ_{k} : x_{1},x_{2},\widetilde{x} \in I\} \\
& = & \sum_{j=0}^{k-3} W_{P_{1} \setminus \{x_{1}\}}\(j\) \cdot 
W_{P_{2} \setminus \{x_{2}\}}\(k-3-j\),
\end{eqnarray*}
and the desired result follows.
\end{proof}

\begin{thm} \label{T:even}
For all $k,v \in \N$ such that $k \leq 2v$, 
\begin{eqnarray*}
f_{2v,k} & = & \sum_{j \geq 0} \sum_{r \geq 0} 
\binom{j}{r} \binom{v-j}{r} \binom{v-\(j+r\)}{k-\(2j+r\)} \\
& = & \sum_{j \geq 0} 
\frac{\(k-2j+1\)_{j}\(v+j-k+1\)_{k-2j}}{j!\(k-2j\)!}
\end{eqnarray*}
holds.
\end{thm}
\begin{proof}
For any $n \in \N \setminus\{0\}$ write the fence poset 
$\cZ_{n}$ as the poset $\{z_{1}, \ldots ,z_{n}\}$ in which 
$z_{2j-1} \triangleleft z_{2j} \triangleright z_{2j+1}$, for all 
$j \geq 1$, are the cover relations, so $\rho\(z_{j}\)=0$ if and 
only if $j \equiv 1 \pmod{2}$ and $\rho\(z_{j}\)=1$ if and only 
if $j \equiv 0 \pmod{2}$. 
If we consider $P_{1}=\cZ_{2v+1}=\{a_{1}, \ldots ,a_{2v+1}\}$, 
$P_{2}=\cZ_{1}=\{b_{1}\}$, we have that 
$P_{1}\(a_{2v+1}\) \circledast P_{2}\(b_{1}\) \simeq \cZ_{2v+3}$, 
and the desired result follows applying Theorems~\ref{T:1fl} 
and~\ref{T:main}, and Chu--Vandermonde summation for hypergeometric 
series as in the proof of Theorem~\ref{T:main}. 
\end{proof}

\begin{cor}
For any $v \in \N$ and all $0 \leq k \leq 2v$, the sequence 
$f_{2v,k}$ is increasing in $v$, viz. $f_{2\(v+1\),k} > f_{2v,k}$. 
\end{cor}
\qed

From Proposition~\ref{P:5} and Theorem~\ref{T:even} we 
immediately get the following result.

\begin{thm} \label{T:cr}
For all $k,n \in \N$ such that $k \leq 2n$, 
\begin{eqnarray*}
c_{n,k} & = & \sum_{j \geq 0} 
\frac{\(k-2j+1\)_{j-2}
\(n+j-k+1\)_{k-2j-2}}{j!\(k-2j\)!} \\
& \cdot & \[\(k-j-1\)_{2}\(n-j-1\)_{2}
-\(\(k-2j-1\)_{2}\)^{2}\]
\end{eqnarray*}
holds.
\end{thm}
\qed

Therefore Theorems~\ref{T:main}, \ref{T:even} and~\ref{T:cr} give the 
solution of the recursive identities in Proposition~\ref{P:5}.

\section{Generalized Fences with higher asymmetric peaks}

\label{S:as}

Now we define an \emph{asymmetric peak} poset with two positive integers 
parameters $\mu,\nu$.

\begin{defn}
Let $\mu,\nu \in \N \setminus\{0\}$; we define the poset 
\emph{asymmetric peak} $\(\operatorname{AP}_{\mu,\nu}, \leq\)$ 
in the following way:
$\operatorname{AP}_{\mu,\nu} = 
\{a_{j} : j \in \[\mu\]\} \biguplus \{b_{j} : j \in \[\nu\]\}
\biguplus \{\omega\}$, and the cover relations are
\begin{description}
\item{$\circ$} $a_{j} \triangleleft a_{j+1}$ for all $j \in \[\mu-1\]$,
\item{$\circ$} $b_{j} \triangleleft b_{j+1}$ for all $j \in \[\nu-1\]$,
\item{$\circ$} $a_{\mu} \triangleleft \omega$,
\item{$\circ$} $b_{\nu} \triangleleft \omega$.
\end{description}
\end{defn}

\begin{prop} \label{P:ap}
Let $\mu,\nu \in \N \setminus\{0\}$; then
$$W_{\operatorname{AP}_{\mu,\nu}}\(k\)= \begin{cases}
1 \qquad & \text{ if } k=0 \text{ or } k=\mu + \nu + 1 \\
k+1 & \text{ if } k \leq \min\{\mu,\nu\} \\
\min\{\mu,\nu\} + 1
    & \text{ if } \min\{\mu,\nu\} \leq k \leq \max\{\mu,\nu\} \\ 
1 + \mu + \nu - k & \text{ if } \max\{\mu,\nu\} \leq k
\end{cases}$$
holds, for any $k=0, \ldots , 
\# \operatorname{AP}_{\mu,\nu} = \mu + \nu + 1$.
\end{prop}
\begin{proof}
The result is clear is $k=0$ or $k=\mu + \nu + 1$. 

We consider the case $\mu \leq \nu$, the case $\mu \geq \nu$ is 
completely symmetric. \\
If $k \in \[\mu+\nu\]$ then any 
$I \in \cJ\(\operatorname{AP}_{\mu,\nu}\)$ with $\rho\(I\)=\# I = k$ 
has the shape 
$I=\{a_{j} : j \in \[r\]\} \biguplus \{b_{j} : j \in \[t\]\}$ with
\begin{equation} \label{E:tr}
r+t=k,
\end{equation}
so $W_{\operatorname{AP}_{\mu,\nu}}\(k\)$ equals
the number of solutions of~\eqref{E:tr} with the constraints
\begin{align*}
0 \leq r \leq k  \; \text{ and } \; 0 \leq t \leq k \qquad
   & \text{ if } k \leq \mu, \\
0 \leq r \leq \mu  \; \text{ and } \; k - \mu \leq t \leq k \qquad
   & \text{ if } \mu \leq k \leq \nu, \\ 
 k - \nu \leq r \leq \mu  \; \text{ and } \; k - \mu \leq t \leq \nu 
   \qquad & \text{ if } \nu \leq k.
\end{align*}

The desired result follows.
\end{proof}

Results proved in \S~\ref{S:cf} allows to get explicit closed formulas 
for Whitney numbers of lattices of order ideals of ``fences with 
higher asymmetric peaks'', i.e. the alternate composition of fences and 
asymmetric peaks by the operator $\circledast$, 
see Definition~\ref{D:star}.

\begin{figure}[ht] 
\begin{center}
\includegraphics[width=12cm]{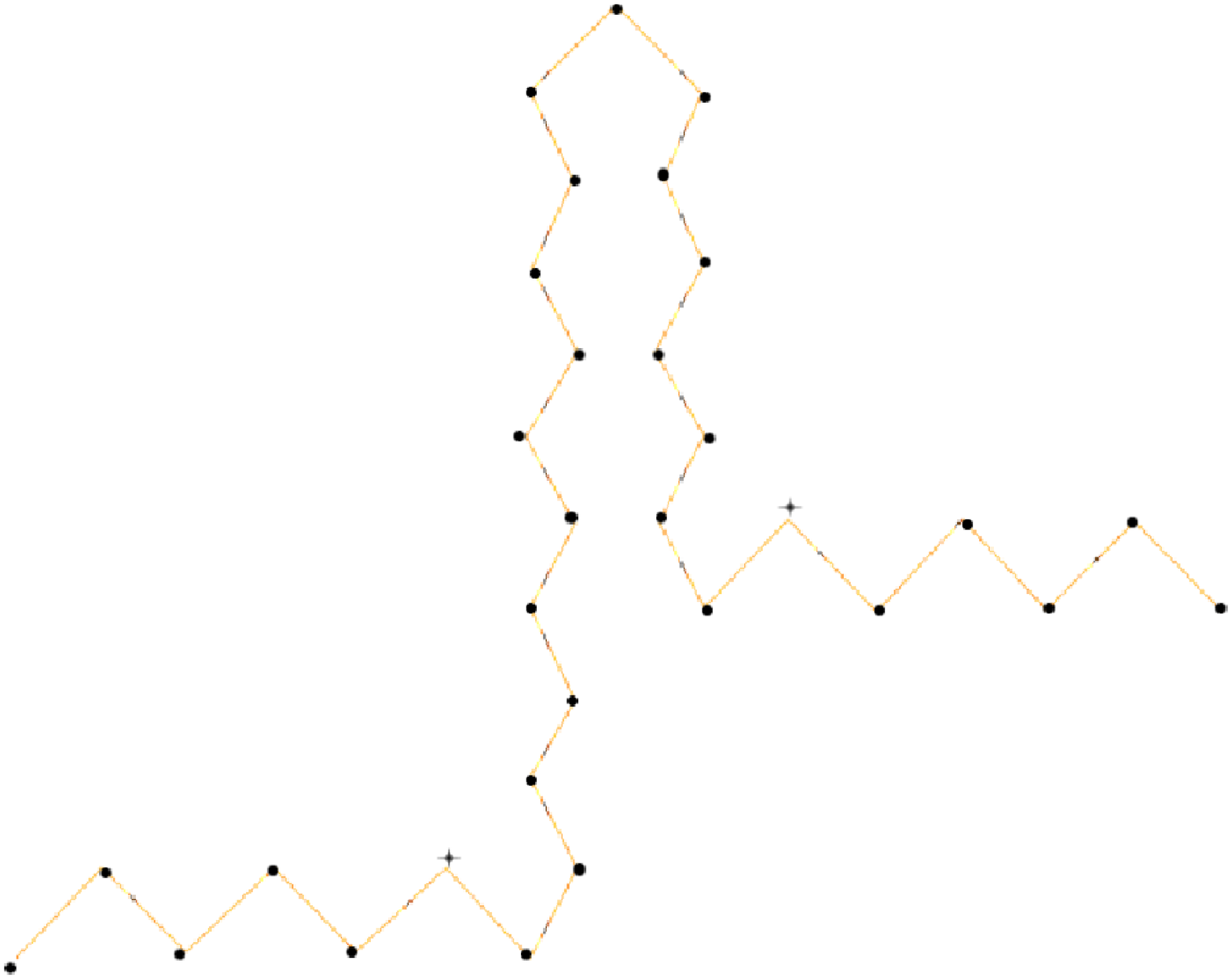}
\caption{$\operatorname{FAP}\(7,10,6,7\)$}
\label{F:as}
\end{center}
\end{figure}

For example, we can consider a fence with one higher asymmetric peak, which 
can be formally defined as the following poset $\(FAP\(w,x,y,z\),\leq\)$ 
where with $w,x,y,z \in \N$ and $w \equiv 1 \pmod{2}$: 
$$\operatorname{FAP}\(w,x,y,z\)=\{a_{1}, \ldots ,a_{w},b_{1}, \ldots ,b_{x},
\omega, c_{1}, \ldots ,c_{y},d_{1}, \ldots ,d_{z}\},$$ 
where the cover relations are 
\begin{description}
\item{$\circ$} $a_{2j-1} \triangleleft a_{2j} \triangleright a_{2j+1}$, 
               for all $j \geq 1$,
\item{$\circ$} $a_{w} \triangleleft b_{1}$,
\item{$\circ$} $b_{j} \triangleleft b_{j+1}$ for all $j \in \[x-1\]$,
\item{$\circ$} $c_{j} \triangleleft c_{j+1}$ for all $j \in \[y-1\]$,
\item{$\circ$} $b_{x} \triangleleft \omega$,
\item{$\circ$} $c_{y} \triangleleft \omega$,
\item{$\circ$} $d_{1} \triangleleft c_{1}$,
\item{$\circ$} $d_{2j-1} \triangleleft d_{2j} \triangleright d_{2j+1}$, 
               for all $j \geq 1$.
\end{description}
In Figure~\ref{F:as} the Hasse diagram of $\operatorname{FAP}\(7,10,6,7\)$ is depicted.

Inside $\operatorname{FAP}\(w,x,y,z\)$ consider the subposets
\begin{eqnarray*}
\cP_{1} & = & \{a_{1}, \ldots ,a_{w-2}\} \simeq \cZ_{w-2}, \\
\cP_{2} & = & \{a_{w},b_{1}, \ldots ,b_{x}, \omega, c_{1}, 
        \ldots ,c_{y}, d_{1}\} \simeq \operatorname{AP}_{x+1,y+1}, \\
\cP_{3} & = & \{d_{3}, \ldots ,d_{z}\} \simeq \cZ_{z-2};
\end{eqnarray*}
therefore 
$\operatorname{FAP}\(w,x,y,z\)=\cP_{1} \biguplus 
\{a_{w-1}\} \biguplus \cP_{2} \biguplus \{d_{2}\} \biguplus \cP_{3}$.

We have that 
\begin{eqnarray*}
\operatorname{FAP}\(w,x,y,z\) & \simeq &
\(\cP_{1}\(a_{w-2}\) \circledast \cP_{2}\(a_{w}\)\)\(d_{1}\)
\circledast \cP_{3}\(d_{3}\) \\ 
& \simeq & \cP_{1}\(a_{w-2}\) \circledast 
\(\cP_{2}\(d_{1}\) \circledast \cP_{3}\(d_{3}\)\)\(a_{w}\),
\end{eqnarray*}
therefore from Theorems~\ref{T:main}, \ref{T:1fl}, 
\ref{T:even}, \ref{T:cr} and Proposition~\ref{P:ap} we get an 
explicit closed formulas for the rank polynomial of the distributive 
lattice of all order ideals of the poset 
$\operatorname{FAP}\(w,x,y,z\)$.

We remark that the same construction can be iterated, so for any 
non--negative integer $k$ we can recursively have a formula for 
the rank polynomial of the lattice of all order ideals of a fence 
with $k$ higher asymmetric peaks.

\section{Open problems and Conjectures}

In~\cite{MZ} using recursive formulas stated in Proposition~\ref{P:5}
it is proved that sequences $f_{n,k}$ and $c_{n,k}$ are indeed unimodal; 
for definitions and comprehensive surveys about unimodal and (strong) 
log--concave sequences we refer to~\cite{Bre2,Hog,Pro,Sta1,Wag} 
(and the references therein).

We feel that the following stronger statement is true.

\begin{conj} \label{C:C}
For any $3 \ne n \in \N \setminus \{0\}$ and all $0 \leq k \leq n$, 
the sequence $f_{n,k}$ is log--concave in $k$, viz. 
$f_{n,k}^{2} \geq f_{n,k-1} f_{n,k+1}$ for any $k \in \[n-1\]$. 

Moreover, for any $4 \leq n \in \N$ and all $0 \leq k \leq 2n$, the 
sequence $c_{n,k}$ is strong log--concave in $k$, viz. 
$c_{n,k}^{2} > c_{n,k-1} c_{n,k+1}$ for any $k \in \[2n-1\]$. 
\end{conj}

Using a computer, Conjecture~\ref{C:C} has been verified for 
distributive lattices of order ideals of fences and crowns, for 
all fences and crowns with cardinality less or equal than $90$.

Moreover, we note that it would be of very great interest to study the 
following much more general problem.

\begin{opb}
Characterize finite posets for which the distributive lattice of 
order ideals is rank (strong) log--concave or just rank unimodal.
\end{opb}

\section*{Aknowledgment}

The author would like to thank Francesco Brenti for suggesting this problem 
and his helpful advice, and Christian Krattenthaler for useful remarks.

\end{document}